\documentclass{article}
\usepackage{amsmath, amsthm, amsfonts, amssymb}
\usepackage{hyperref}

\newtheorem{theorem}{Theorem}[section]

\newtheorem{proposition}[theorem]{Proposition}
\newtheorem{corollary}[theorem]{Corollary}
\theoremstyle{definition}
\newtheorem{definition}[theorem]{Definition}

\theoremstyle{remark}

\numberwithin{equation}{section}

\newcommand{\ket}[1]{ | #1  \rangle}

\newcommand{\ii}{\mathbf{i_1}}
\newcommand{\iii}{\mathbf{i_2}}
\newcommand{\jj}{\mathbf{j}}
\newcommand{\ee}{\mathbf{e_1}}
\newcommand{\eee}{\mathbf{e_2}}
\newcommand{\e}[1]{\mathbf{e_{#1}}}
\newcommand{\ik}[1]{\mathbf{i_{#1}}}
\newcommand{\T}{\mathbb{T}}
\newcommand{\C}{\mathbb{C}}
\newcommand{\D}{\mathbb{D}}
\newcommand{\R}{\mathbb{R}}
\newcommand{\nc}{\mathcal{NC}}
\renewcommand{\(}{\left(}
\renewcommand{\)}{\right)}
\newcommand{\oa}{\left\{}
\newcommand{\fa}{\right\}}

\renewcommand{\P}[2]{P_{#1}( #2 )}

\renewcommand{\exp}[1]{\textrm{exp}\oa #1 \fa}

\newcommand{\bo} {\ensuremath{{\bf i_1}}}

\newcommand{\eo} {\ensuremath{{\bf e_1}}}
\newcommand{\et} {\ensuremath{{\bf e_2}}}
\newcommand{\iic}{{\rm \bf{i}}_{\bf{1}}^2}
\newcommand{\iiic}{{\rm \bf{i}}_{\bf{2}}^2}
\newcommand{\h}[1]{{\widehat{#1}}}
\newcommand{\hh}{{\widehat{1}}}
\newcommand{\hhh}{{\widehat{2}}}
\newcommand{\mC}{\ensuremath{\mathbb{C}}}

\newcommand{\rmd}{\mbox{d}}
\newcommand{\rme}{\mbox{e}}
\newcommand{\mo}{\mathbf{1}}
\newcommand{\mt}{\mathbf{2}}
\newcommand{\mk}{\mathbf{k}}
\newcommand{\scalarmath}[2]{\( #1, #2 \)}

\title{\textbf{Infinite-Dimensional Bicomplex Hilbert Spaces}}

\author{\textbf{Rapha\"el Gervais Lavoie$^1$,
Louis Marchildon$^1$} \\ \textbf{and Dominic Rochon$^2$}}

\begin{document}
\date{}
\maketitle

\begin{center}
$^{1}$ D\'epartement de physique, Universit\'e du Qu\'ebec,\\ Trois-Rivi\`eres, Qc. Canada G9A 5H7

\bigskip
$^{2}$ D\'epartement de math\'ematiques et d'informatique,
Universit\'e du Qu\'ebec,\\ Trois-Rivi\`eres, Qc. Canada G9A 5H7

\bigskip
email: raphael.gervaislavoie@uqtr.ca, louis.marchildon@uqtr.ca, dominic.rochon@uqtr.ca.\\
~
\end{center}

\begin{abstract}
This paper begins the study of infinite-dimensional
modules defined on bicomplex numbers.  It generalizes
a number of results obtained with finite-dimensional
bicomplex modules.  The central concept introduced
is the one of a bicomplex Hilbert space.  Properties
of such spaces are obtained through properties of
several of their subsets which have the structure of
genuine Hilbert spaces.  In particular, we derive the Riesz
representation theorem for bicomplex continuous linear
functionals and a general version of the bicomplex Schwarz
inequality.  Applications to concepts relevant to quantum
mechanics, specifically the bicomplex analogue of the quantum
harmonic oscillator, are pointed out.
\end{abstract}
\section{Introduction}

The mathematical structure of quantum mechanics
consists in Hilbert spaces defined over the field
of complex numbers~\cite{Neumann}.  This structure
has been extremely successful in explaining vast
amounts of experimental data pertaining largely,
but not exclusively, to the world of molecular,
atomic and subatomic phenomena.

Bicomplex numbers~\cite{Price}, just like quaternions,
are a generalization of complex numbers by means of
entities specified by four real numbers.  These two
number systems, however, are different in two
important ways: quaternions, which form a division
algebra, are noncommutative, whereas bicomplex numbers
are commutative but do not form a division algebra.

Division algebras do not have zero divisors, that is,
nonzero elements whose product is zero.  Many believe
that any attempt to generalize quantum mechanics
to number systems other than complex numbers should
retain the division algebra property.  Indeed
considerable work has been done over the years on
quaternionic quantum mechanics~\cite{Adler}.
However, in the past few years it was pointed out that several
features of quantum mechanics can be generalized to bicomplex
numbers.  A generalization of Schr\"{o}dinger's equation
for a particle in one dimension was proposed~\cite{Rochon2}
and self-adjoint operators were defined on finite-dimensional
bicomplex Hilbert spaces~\cite{GMR2, Rochon3}.

In this spirit, eigenvalues and eigenkets or eigenfunctions
of the bicomplex analogue of the quantum harmonic oscillator
Hamiltonian were obtained in full generality~\cite{GMR}
over an infinite-dimensional bicomplex module~$M$.
The harmonic oscillator is one of the simplest and, at the same
time, one of the most important systems of quantum mechanics,
involving as it is an infinite-dimensional vector space.
However, the module~$M$ defined in~\cite{GMR} does not have a
property of completeness. Indeed it is made up of finite linear
combinations of eigenkets or, in the coordinate basis, of
hyperbolic Hermite polynomials.  Hence Cauchy sequences of
elements of~$M$ do not in general converge to an element
of~$M$.

In this paper we introduce the mathematical tools
necessary to investigate the bicomplex analogue
of quantum-mechanical state spaces.  This we do by
defining the concept of infinite-dimensional bicomplex Hilbert
space, closely related to its complex version
and already introduced in~\cite{GMR2, Rochon3} for finite
dimensions.  In Section~2, we summarize known algebraic
properties of bicomplex numbers and recall the concepts
of bicomplex modules and scalar products.  Section~3
contains the main results of this paper.  Bicomplex
Hilbert spaces are introduced and subsets are identified
that have the structure of standard complex Hilbert
spaces.  We derive a general version of the bicomplex
Schwarz inequality and a version of Riesz's representation
theorem for bicomplex continuous linear functionals.
For an arbitrary bicomplex Hilbert space $M$, the dual space
$M^{*}$ of continuous linear functionals on $M$ can then be
identified with~$M$ through the bicomplex scalar product
$\scalarmath{\cdot}{\cdot}$.  These are the tools necessary
to better justify the concept of self-adjoint operators
acting in an infinite-dimensional bicomplex Hilbert space
introduced in~\cite{GMR}.  Bicomplex Hilbert spaces with
countable bases are discussed and results proved on
important Hilbert subspaces.  Section~4 examines the
example of the quantum harmonic oscillator.  In the standard
case the Hamiltonian eigenfunctions generate the state
space $V=L^2(\mathbb{R})$.  We construct an infinite-dimensional
bicomplex module~$M=(\eo V)\oplus (\et V)$, which has the property
of completeness and has the eigenfunctions of the
bicomplex harmonic oscillator Hamiltonian as a basis.

\section{Preliminaries}\label{Preliminaries}

This section first summarizes a number of known results
on the algebra of bicomplex numbers, which will be needed
in this paper.  Much more details as well as proofs can
be found in~\cite{Price, Rochon1, Rochon2, Rochon3}.
Basic definitions related to bicomplex modules and scalar
products are also formulated as in~\cite{GMR2, Rochon3}, but
here we make no restrictions to finite dimensions.


\subsection{Bicomplex Numbers}\label{Bicomplex Numbers}

\subsubsection{Definition}\label{Definition of bicomplex numbers}
The set $\T$ of \emph{bicomplex numbers} is defined as
\begin{align}
\mathbb{T}:=\{ w=z_1+z_2\mathbf{i_2}~|~z_1,z_2\in\mathbb{C}(\mathbf{i_1}) \},
\label{2.1}
\end{align}
where $\ii$ and $\iii$ are independent imaginary units such that
$\iic=-1=\iiic$.  The product of $\ii$ and $\iii$ defines
a hyperbolic unit $\jj$ such that $\mathbf{j}^2=1$.
The product of all units is commutative and satisfies
\begin{equation}
\ii\iii=\jj, \qquad \ii\jj=-\iii,
\qquad \iii\jj=-\ii. \label{2.2}
\end{equation}
With the addition and multiplication of two
bicomplex numbers defined in the obvious way,
the set $\mathbb{T}$ makes up a commutative ring.

Three important subsets of $\mathbb{T}$ can be
specified as
\begin{align}
\mathbb{C}(\ik{k}) &:= \{ x+y\ik{k}~|~x,y\in\mathbb{R} \},
\qquad k=1,2 ;\label{2.3}\\
\mathbb{D} &:= \{ x+y\jj~|~x,y\in\mathbb{R} \} .\label{2.4}
\end{align}
Each of the sets $\mathbb{C}(\ik{k})$ is isomorphic
to the field of complex numbers, while $\mathbb{D}$ is
the set of so-called \emph{hyperbolic numbers}.

\subsubsection{Conjugation and Moduli}\label{Bicomplex conjugation}

Three kinds of conjugation can be defined on
bicomplex numbers. With $w$ specified as in~\eqref{2.1}
and the bar ($\,\bar{\mbox{}}\,$) denoting complex
conjugation in $\mathbb{C}(\mathbf{i_1})$,
we define
\begin{equation}
w^{\dag_1}:=\bar{z}_1+\bar{z}_2\mathbf{i_2},\label{2.5}
\qquad w^{\dag_2}:=z_1-z_2\mathbf{i_2},
\qquad w^{\dag_3}:=\bar{z}_1-\bar{z}_2\mathbf{i_2} .
\end{equation}
It is easy to check that each conjugation has the following
properties:
\begin{equation}
(s+t)^{\dag_k}=s^{\dag_k}+t^{\dag_k},
\qquad \left(s^{\dag_k} \right)^{\dag_k}=s,
\qquad (s\cdot t)^{\dag_k}=s^{\dag_k}\cdot t^{\dag_k} .
\label{2.6}
\end{equation}
Here $s,t\in\mathbb{T}$ and $k=1,2,3$.

With each kind of conjugation, one can define a specific
bicomplex modulus as
\begin{align}
|w|_\ii^2&:=w\cdot w^{\dag_2}=z_1^2+z_2^2~\in\C(\ii),\label{2.7a}\\
|w|_\iii^2&:=w\cdot w^{\dag_1}=\left(|z_1|^2-|z_2|^2\right)
+ 2 \, \textrm{Re}(z_1\bar{z}_2)\iii~\in\C(\iii),\label{2.7b}\\
|w|_\jj^2&:=w\cdot w^{\dag_3}=\left(|z_1|^2+|z_2|^2\right)
- 2 \, \textrm{Im}(z_1\bar{z}_2)\jj~\in\D.\label{2.7c}
\end{align}
It can be shown that $|s\cdot t|_k^2=|s|_k^2\cdot|t|_k^2$,
where $k=\ii,\iii$ or $\jj$.

In this paper we will often use the Euclidean $\R^4$ norm
defined as
\begin{equation}
|w|:=\sqrt{|z_1|^2+|z_2|^2}=\sqrt{\textrm{Re}(|w|_\jj^2)} \; .
\label{2.8}
\end{equation}
Clearly, this norm maps $\T$ into $\R$.  We have $|w|\geq0$,
and $|w|=0$ if and only if $w=0$. Moreover~\cite{Rochon1},
for all $s,t\in\T$,
\begin{equation}
|s+t|\leq|s|+|t|, \qquad |s\cdot t|\leq \sqrt{2} \, |s|\cdot|t|.
\label{2.9}
\end{equation}

\subsubsection{Idempotent Basis}\label{Idempotant basis}

Bicomplex algebra is considerably simplified by
the introduction of two bicomplex numbers $\ee$
and $\eee$ defined as
\begin{equation}
\ee:=\frac{1+\jj}{2},\qquad\eee:=\frac{1-\jj}{2}.\label{2.10}
\end{equation}
In fact $\ee$ and $\eee$ are hyperbolic numbers.
They make up the so-called \emph{idempotent basis}
of the bicomplex numbers. One easily checks that ($k=1,2$)
\begin{equation}
\mathbf{e}_{\mathbf{1}}^2=\ee,
\quad \mathbf{e}_{\mathbf{2}}^2=\eee,
\quad \ee+\eee=1,
\quad \mathbf{e}_{\mathbf{k}}^{\dag_3}=\e{k} ,
\quad \ee\eee=0 . \label{2.11}
\end{equation}

Any bicomplex number $w$ can be written uniquely as
\begin{equation}
w = z_1+z_2\iii = z_\hh \ee + z_\hhh \eee , \label{2.12}
\end{equation}
where
\begin{equation}
z_\hh= z_1-z_2\ii \quad \mbox{and}
\quad z_\hhh= z_1+z_2\ii \label{2.12a}
\end{equation}
both belong to $\mathbb{C}(\ii)$.  Note that
\begin{equation}
|w| = \frac{1}{\sqrt{2}}
\sqrt{|z_\hh |^2 + |z_\hhh |^2} \, . \label{norm7}
\end{equation}
The caret notation
($\hh$ and $\hhh$) will be used systematically in
connection with idempotent decompositions, with the
purpose of easily distinguishing different types
of indices.  As a consequence of~\eqref{2.11}
and~\eqref{2.12}, one can check that if
$\sqrt[n]{z_\hh}$ is an $n$th root of $z_\hh$
and $\sqrt[n]{z_\hhh}$ is an $n$th root of $z_\hhh$,
then $\sqrt[n]{z_\hh} \, \ee + \sqrt[n]{z_\hhh} \, \eee$
is an $n$th root of $w$.

The uniqueness of the idempotent decomposition
allows the introduction of two projection operators as
\begin{align}
P_1: w \in\T&\mapsto z_\hh \in\C(\ii),\label{2.14}\\
P_2: w \in\T&\mapsto z_\hhh \in\C(\ii).\label{2.15}
\end{align}
The $P_k$ ($k = 1, 2$) satisfy
\begin{equation}
[P_k]^2=P_k, \qquad P_1\ee+P_2\eee=\mathbf{Id}, \label{2.16}
\end{equation}
and, for $s,t\in\T$,
\begin{equation}
P_k(s+t)=P_k(s)+P_k(t),
\qquad P_k(s\cdot t)=P_k(s)\cdot P_k(t) .\label{2.17}
\end{equation}

The product of two bicomplex numbers $w$ and $w'$
can be written in the idempotent basis as
\begin{align}
w \cdot w' = (z_\hh \ee + z_\hhh \eee)
\cdot (z'_\hh \ee + z'_\hhh \eee)
= z_\hh z'_\hh \ee + z_\hhh z'_\hhh \eee .\label{2.20}
\end{align}
Since 1 is uniquely decomposed as $\ee + \eee$,
we can see that $w \cdot w' = 1$ if and only if
$z_\hh z'_\hh = 1 = z_\hhh z'_\hhh$.  Thus $w$ has an inverse
if and only if $z_\hh \neq 0 \neq z_\hhh$, and the
inverse $w^{-1}$ is then equal to
$(z_\hh)^{-1} \ee + (z_\hhh)^{-1} \eee$.  A nonzero $w$ that
does not have an inverse has the property that
either $z_\hh = 0$ or $z_\hhh = 0$, and such a $w$ is
a divisor of zero.  Zero divisors make up the
so-called \emph{null cone} $\nc$.  That terminology comes
from the fact that when $w$ is written as in~\eqref{2.1},
zero divisors are such that $z_1^2 + z_2^2 = 0$.

Any hyperbolic number can be written in the
idempotent basis as $x_\hh \ee + x_\hhh \eee$, with
$x_\hh$ and $x_\hhh$ in~$\R$.  We define the set~$\D^+$
of positive hyperbolic numbers as
\begin{equation}
\D^+:= \{ x_\hh \ee + x_\hhh \eee ~|~ x_\hh, x_\hhh \geq 0 \}.
\label{2.21}
\end{equation}
Since $w^{\dag_3} = \bar{z}_\hh \ee + \bar{z}_\hhh \eee$,
it is clear that $w \cdot w^{\dag_3} \in \D^+$ for any
$w$ in $\T$.

\subsection{$\T$-Modules}\label{Module}

The set of bicomplex numbers is a commutative ring.
Just like vector spaces are defined over fields,
modules are defined over rings \cite{Bourbaki}.  A module~$M$ defined
over the ring $\T$ of bicomplex numbers is called a
$\T$-\emph{module}~\cite{GMR2, Rochon3}.

\begin{definition}
Let $M$ be a $\T$-module. For $k=1, 2$, we define $V_k$
as the set of all elements of the form $\e{k} \ket{\psi}$,
with $\ket{\psi} \in M$.  Succinctly, $V_1:=\eo M$
and $V_2:=\et M$.
\label{Definition3.1}
\end{definition}

\noindent
For $k=1,2$, addition and multiplication by a $\mC(\bo)$ scalar
are closed in $V_k$.  Therefore, $V_k$ is a vector space over
$\mC(\bo)$.  Any element $\ket{v_k} \in V_k$ satisfies
$\ket{v_k} = \e{k} \ket{v_k}$.

For arbitrary $\T$-modules, vector spaces $V_1$ and $V_2$
bear no structural similarities.  For more specific modules,
however, they may share structure.  It was shown in~\cite{GMR2}
that if~$M$ is a finite-dimensional free $\T$-module, then
$V_1$ and~$V_2$ have the same dimension.  Other similar
instances will be examined in Section~3.

\begin{proposition}
Let $M$ be a $\T$-module and let $\ket{\psi}\in M$.
There exist unique vectors $\ket{v_1}\in V_1$
and $\ket{v_2}\in V_2$ such that $\ket{\psi}=\ket{v_1}+\ket{v_2}$.
\label{decomp}
\label{Proposition3.1}
\end{proposition}
\begin{proof}
Let $\ket{\psi}\in M$.  We can always write
$$\ket{\psi}=\eo\ket{\psi}+\et\ket{\psi}=\ket{v_1}+\ket{v_2},$$
where $\ket{v_k}:=\e{k}\ket{\psi}\in V_k$, for $k=1,2$.
Suppose that $\ket{\psi}=\ket{v'_1}+\ket{v'_2}$,
with $\ket{v'_k}\in V_k$. Then
$$\ket{v_1}+\ket{v_2}=\ket{v'_1}+\ket{v'_2} .$$
Multiplying both sides with $\e{k}$ and making use
of \eqref{2.11}, we obtain
$$\ket{v_k} = \e{k}\ket{v_k} = \e{k}\ket{v'_k} = \ket{v'_k}$$
for $k=1,2$.
\end{proof}

\medskip\noindent
Henceforth we will write $\ket{\psi}_\mk = \e{k} \ket{\psi}$,
keeping in mind that $\e{k} \ket{\psi}_\mk = \ket{\psi}_\mk$.
Proposition~\ref{Proposition3.1} immediately leads to
the following result.
\begin{theorem}
The $\T$-module $M$ can be viewed as a vector space $M'$
over $\mC(\bo)$, and $M'=V_1\oplus V_2.$
\label{SD}
\label{Theo3.3}
\end{theorem}
From a set-theoretical point of view, $M$ and $M'$ are
identical.  In this sense we can say, perhaps improperly,
that the \textbf{module} $M$ can be decomposed into the
direct sum of two vector spaces over $\mC(\bo)$, i.e.\
$M=V_1\oplus V_2.$

\subsection{Bicomplex Scalar Product}

The norm of a vector is an important concept in
vector space theory.  We will now generalize it to
$\T$-modules, making use of the association
established in Theorem~\ref{Theo3.3}.

\begin{definition}
Let $M$ be a $\T$-module and let $M'$ be the associated
vector space. We say that $\|\cdot\|:M\longrightarrow \mathbb{R}$
is a \textbf{$\T$-norm} on $M$ if the following holds:

\smallskip\noindent
1. $\|\cdot\|:M'\longrightarrow \mathbb{R}$ is a norm;\\
2. $\big{\|}w\cdot \ket{\psi}\big{\|}\leq \sqrt{2}
\big{|}w\big{|}\cdot\big{\|}\ket{\psi}\big{\|}$,
$\forall w\in\T$, $\forall \ket{\psi}\in M$.
\label{norm}
\end{definition}
\noindent A $\T$-module with a \textbf{$\T$-norm} is called a
\textbf{normed $\T$-module}.

In vector space theory, a norm can be induced by a scalar
product.  Having in mind the use of such norms, we recall
the definition of a \textbf{bicomplex scalar product}
introduced in~\cite{Rochon3} (the physicists' ordering
convention being used).

\begin{definition}
Let $M$ be a $\mathbb{T}$-module.  Suppose that with
each pair $\ket{\psi}$ and $\ket{\phi}$ in $M$, taken in this
order, we associate a bicomplex number $\(\ket{\psi},\ket{\phi}\)$.
We say that the association defines a bicomplex scalar (or
inner) product if it satisfies the following conditions:

\smallskip\noindent
1. $(\ket{\psi},\ket{\phi}+\ket{\chi})
=(\ket{\psi},\ket{\phi})+(\ket{\psi},\ket{\chi})$,
$\forall \ket{\psi},\ket{\phi},\ket{\chi}\in M;$\\
2. $(\ket{\psi},\alpha \ket{\phi})=\alpha (\ket{\psi},\ket{\phi}),$
$\forall \alpha\in\mathbb{T}$, $\forall \ket{\psi},\ket{\phi}\in M;$ \\
3. $(\ket{\psi},\ket{\phi})=(\ket{\phi},\ket{\psi})^{\dagger_{3}}$,
$\forall \ket{\psi},\ket{\phi}\in M;$\\
4. $(\ket{\psi},\ket{\psi})=0 \mbox{ }
\Leftrightarrow \mbox{ }\ket{\psi}=0$, $\forall \ket{\psi}\in M.$
\label{scalar}
\label{Definition4.1}
\end{definition}
%
Property~3 implies that $\(\ket{\psi},\ket{\psi}\)\in\mathbb{D}$.
Definition~\ref{scalar} is intended to be
very general.  In this paper we shall
be more restrictive, by requiring the bicomplex scalar product
$\(\cdot,\cdot\)$ to be \textit{hyperbolic positive}, that is,
\begin{equation}
(\ket{\psi},\ket{\psi})\in\mathbb{D}^{+}, \;
\forall\ket{\psi}\in M. \label{hyperpositive}
\end{equation}

\noindent From Definition~\ref{Definition4.1} it is easy to see that the following projection of a bicomplex scalar product:
\begin{equation}
(\cdot,\cdot)_{\widehat{k}}:=P_k((\cdot,\cdot)):M\times M\longrightarrow \mC(\bo)
\end{equation}
is a \textbf{standard scalar product} on $V_k$, for $k=1,2$.

\begin{theorem}
Let $\ket{\psi},\ket{\phi}\in M$, then
\begin{equation*}
(\ket{\psi},\ket{\phi})
=\eo(\ket{\psi}_\mo,\ket{\phi}_\mo)_{\widehat{1}}
+ \et(\ket{\psi}_\mt,\ket{\phi}_\mt)_{\widehat{2}}
\label{decompEq}
\end{equation*}
\label{decomp2}\label{Theo4.1}
\end{theorem}

\vspace{-5ex}\noindent
\begin{proof}
Since $\ket{\psi}=\ket{\psi}_\mo + \ket{\psi}_\mt$ and
$\ket{\phi}= \ket{\phi}_\mo+\ket{\phi}_\mt$, we have
\begin{align*}
\scalarmath{\ket{\psi}}{\ket{\phi}}
&=\scalarmath{\ket{\psi}_\mo + \ket{\psi}_\mt}{\ket{\phi}_\mo + \ket{\phi}_\mt}\\
&=\eo\scalarmath{\ket{\psi}_\mo} {\ket{\phi}_\mo}
+ \et\scalarmath{\ket{\psi}_\mt}{\ket{\phi}_\mt}\\
&=\ee\oa \ee\P{1}{\scalarmath{\ket{\psi}_\mo}{\ket{\phi}_\mo}}+\eee\P{2}{\scalarmath{\ket{\psi}_\mo}{\ket{\phi}_\mo}} \fa\\
&\qquad+\eee\oa \ee\P{1}{\scalarmath{\ket{\psi}_\mt}{\ket{\phi}_\mt}}
+\eee\P{2}{\scalarmath{\ket{\psi}_\mt}{\ket{\phi}_\mt}} \fa\\
&=\ee\P{1}{\scalarmath{\ket{\psi}_\mo}{\ket{\phi}_\mo}}
+ \eee\P{2}{\scalarmath{\ket{\psi}_\mt}{\ket{\phi}_\mt}}\\
&=\ee\scalarmath{\ket{\psi}_\mo}{\ket{\phi}_\mo}_\hh
+\eee\scalarmath{\ket{\psi}_\mt}{\ket{\phi}_\mt}_\hhh.
\end{align*}
\end{proof}

\noindent
We point out that a bicomplex scalar product is
\textbf{completely characterized} by the two standard
scalar products $\scalarmath{\cdot}{\cdot}_{\widehat{k}}$ on $V_k$.
In fact, if $\scalarmath{\cdot}{\cdot}_{\widehat{k}}$
is an arbitrary scalar product on $V_k$, for $k=1,2$,
then $\scalarmath{\cdot}{\cdot}$ defined as in \eqref{decompEq}
is a bicomplex scalar product on $M$.

\section{Bicomplex Hilbert Spaces}

\subsection{General Results}

In this section we define the notion of a bicomplex
Hilbert space and prove the analog of the Riesz
representation theorem.

\begin{definition}
Let $M$ be a $\mathbb{T}$-module and let $(\cdot,\cdot)$
be a bicomplex scalar product defined on $M$. The space
$\{M, (\cdot,\cdot)\}$ is called a $\mathbb{T}$-inner product
space, or bicomplex pre-Hilbert space.  When no confusion
arises, $\{M, (\cdot,\cdot)\}$ will simply be denoted by~$M$.
\end{definition}

\begin{theorem}
Let $M$ be a bicomplex pre-Hilbert space. Then
$(V_k, \scalarmath{\cdot}{\cdot}_{\widehat{k}})$ is a complex
$($in $\mC(\bo))$ pre-Hilbert space for $k=1,2$.
\label{Pre}
\label{Theo4.2}
\end{theorem}
\begin{proof}
Since $(\cdot,\cdot)_{\widehat{k}}$ is a standard scalar product
when $M'$ is restricted to the vector space $V_k$, then
$(V_k, \scalarmath{\cdot}{\cdot}_{\widehat{k}})$ is a complex
(in $\mC(\bo)$) pre-Hilbert space.
\end{proof}


\medskip\noindent
If $V_1$ and $V_2$ are complete, then $M'=V_1\oplus V_2$
is a direct sum of two Hilbert spaces.  It is easy to see that
$M'$ is also a Hilbert space, when the following natural
scalar product is defined over the direct sum~\cite{Conway}:
\begin{equation}
\scalarmath{\ket{\psi}_\mo \oplus \ket{\psi}_\mt}
{\ket{\phi}_\mo \oplus \ket{\phi}_\mt}
= \scalarmath{\ket{\psi}_\mo}{\ket{\phi}_\mo}_{\widehat{1}}
+ \scalarmath{\ket{\psi}_\mt}{\ket{\phi}_\mt}_{\widehat{2}} .
\label{4.6}
\end{equation}
From this scalar product, we can define a \textbf{norm}
on the vector space $M'$:
\begin{align}
\big{|}\big{|}\ket{\phi}\big{|}\big{|}
&:= \frac{1}{\sqrt{2}}
\sqrt{\scalarmath{\ket{\phi}_\mo} {\ket{\phi}_\mo}_{\widehat{1}}
+ \scalarmath{\ket{\phi}_\mt}{\ket{\phi}_\mt}_{\widehat{2}}} \notag\\
&=\frac{1}{\sqrt{2}} \sqrt{ \big{|}\ket{\phi}_\mo\big{|}^{2}_{1}
+ \big{|}\ket{\phi}_\mt\big{|}^{2}_{2}} \, .
\label{T-norm}
\end{align}
Here we wrote
\begin{equation}
\big{|}\ket{\phi}_\mk \big{|}_{k}
= \sqrt{\scalarmath{\ket{\phi}_\mk}
{\ket{\phi}_\mk}_{\widehat{k}}} \, ,
\label{normk1}
\end{equation}
where $|\cdot|_k$ is the natural scalar-product-induced norm on~$V_k$.
The $1/\sqrt{2}$ factor in \eqref{T-norm} is introduced so as to
relate in a simple manner the norm with the bicomplex scalar product.
Indeed we have
\begin{equation}
\big{|}\big{|}\ket{\phi}\big{|}\big{|}
= \frac{1}{\sqrt{2}}
\sqrt{\scalarmath{\ket{\phi}_\mo} {\ket{\phi}_\mo}_{\widehat{1}}
+ \scalarmath{\ket{\phi}_\mt}{\ket{\phi}_\mt}_{\widehat{2}}}
= \big{|} \sqrt{\scalarmath{\ket{\phi}}{\ket{\phi}}} \big{|} ,
\label{T-norma}
\end{equation}
which is easily seen through \eqref{decompEq},
\eqref{norm7} and the remark on roots made after that
last equation.

It is easy to check that $\|\cdot\|$ is a \textbf{$\T$-norm} on $M$
and that the $\T$-module $M$ is \textbf{complete} with respect
to the following metric on $M$:
\begin{equation}
d(\ket{\phi},\ket{\psi})=\big{|}\big{|}\ket{\phi}-\ket{\psi}\big{|}\big{|}.
\end{equation}
Thus $M$ is a \textbf{complete $\T$-module}.

Let us summarize what we found by means of a definition
and a theorem.

\begin{definition}
A bicomplex Hilbert space is a $\T$-inner product space $M$ which is complete with respect to the induced $\T$-norm \eqref{T-norm}.
\label{Hilbert}
\end{definition}

\begin{theorem}
Let $\{M, (\cdot,\cdot)\}$ be a bicomplex
pre-Hilbert space, and let the induced
space $V_k$ be complete with respect to the inner product
$\scalarmath{\cdot}{\cdot}_{\widehat{k}}$ for $k=1, 2$.
Then $\{M, (\cdot,\cdot)\}$ is a bicomplex Hilbert space.
\label{HilbertT}
\end{theorem}

\noindent
The converse of this theorem is also true, as shown by the
following two results.

\begin{theorem}
Let $M$ be a bicomplex Hilbert space.
The normed vector space $(V_k,||\cdot||)$ is closed in $M$ for $k=1,2$.
\label{closed}
\end{theorem}
\begin{proof}
Let $\{\ket{\psi_n}\}\in V_k$ be a sequence converging to a ket
$\ket{\psi}\in M$. Using Property~2 of Definition~\ref{norm},
with $w = \e{k}$,
we see that $\ket{\psi_n}_\mk \rightarrow \ket{\psi}_\mk$
when $n \rightarrow \infty$.  Thus $\ket{\psi}_\mk = \ket{\psi}$
and the sequence $\{\ket{\psi_n}\}$ converges in $V_k$.
Hence $(V_k,||\cdot||)$ is closed in $M$.
\end{proof}

\begin{corollary}
Let $M$ be a bicomplex Hilbert space. Then
$(V_k, \scalarmath{\cdot}{\cdot}_{\widehat{k}})$ is a
complex $($in $\mC(\bo))$ Hilbert space for $k=1,2$.
\label{Hilbert1}
\label{Corollary4.1}
\end{corollary}
\begin{proof}
By Theorem \ref{Pre}, $(V_k, \scalarmath{\cdot}{\cdot}_{\widehat{k}})$
is a normed space over $\mC(\bo)$.  Using the $\T$-norm introduced
in \eqref{T-norm} and Theorem~\ref{closed}, it is easy to see that
$(V_k, |\cdot|_{k})$ is complete since
\begin{equation}
\big{|} \ket{\phi}_\mk \big{|}_{k}
= \sqrt{2} \big{\|} \ket{\phi}_\mk \big{\|} , \;\;
\forall \ket{\phi}_\mk \in V_k . \notag
\end{equation}
Hence $V_k$ is a complex (in $\mC(\bo)$) Hilbert space.
\end{proof}

\medskip\noindent
As a direct application of this result, we obtain the following
\textbf{bicomplex Riesz representation theorem}.

\begin{theorem}[Riesz]
Let $\{M, (\cdot,\cdot)\}$ be a bicomplex Hilbert space
and let $f:M\rightarrow \mathbb{T}$ be a continuous linear
functional on $M$.  Then there is a unique $\ket{\psi}\in M$
such that $\forall \ket{\phi} \in M$,
$f(\ket{\phi})=\scalarmath{\ket{\psi}}{\ket{\phi}}$.
\label{Riesz}
\end{theorem}
\begin{proof}
From the functional $f$ on $M$ we can define a
functional $P_k(f)$ on $V_k$ ($k=1, 2$) as
\begin{equation}
P_k(f)(\ket{\phi}_\mk) := P_k(f(\ket{\phi}_\mk)) . \notag
\end{equation}
It is clear that $P_k(f)$ is linear and continuous
if $f$ is.

We now apply Riesz's theorem~\cite{Herget} to $P_k(f)$
and find a unique $\ket{\psi_k}\in V_k$ such that
$\forall \ket{\phi}_\mk \in V_k$,
\begin{equation}
P_k(f)(\ket{\phi}_\mk) = \scalarmath{\ket{\psi_k}}
{\ket{\phi}_\mk}_{\widehat{k}}. \notag
\end{equation}
Letting $\ket{\psi}:=\ket{\psi_1} + \ket{\psi_2}$
and making use of Theorem~\ref{decomp2}, we get
for every $\ket{\phi}$ in~$M$:
\begin{align*}
(\ket{\psi},\ket{\phi})
&= \eo(\ket{\psi_1}, \ket{\phi}_\mo)_{\widehat{1}}
+ \et(\ket{\psi_2}, \ket{\phi}_\mt)_{\widehat{2}} \\
&= \eo P_1(f) (\ket{\phi}_\mo)
+ \et P_2(f) (\ket{\phi}_\mt) \\
&= \eo P_1 (f(\ket{\phi}_\mo))
+ \et P_2 (f(\ket{\phi}_\mt)) .
\end{align*}
On the other hand,
\begin{align*}
f(\ket{\phi})
&= \eo P_1 (f(\ket{\phi})) + \et P_2 (f(\ket{\phi})) \\
&= \eo P_1 (f( \eo \ket{\phi}_\mo + \et \ket{\phi}_\mt ))
+ \et P_2 (f( \eo \ket{\phi}_\mo + \et \ket{\phi}_\mt )) \\
&= \eo P_1 ( \eo f(\ket{\phi}_\mo) + \et f(\ket{\phi}_\mt ))
+ \et P_2 ( \eo f(\ket{\phi}_\mo) + \et f(\ket{\phi}_\mt )) \\
&= \eo P_1 (f(\ket{\phi}_\mo))
+ \et P_2 (f(\ket{\phi}_\mt)) .
\end{align*}
Hence $f(\ket{\phi}) = (\ket{\psi}, \ket{\phi})$.  If $\ket{\psi '}$
works also, then $(\ket{\psi}, \ket{\phi}) = (\ket{\psi '}, \ket{\phi})$
for any $\ket{\phi}$ in~$M$, and therefore
$(\ket{\psi} - \ket{\psi '}, \ket{\phi}) = 0$.
Letting $\ket{\phi} = \ket{\psi} - \ket{\psi '}$,
we conclude that $\ket{\psi} - \ket{\psi '} = 0$,
that is, $\ket{\psi}$ is unique.
\end{proof}

\medskip\noindent
Theorem~\ref{Riesz} means that for an arbitrary bicomplex
Hilbert space $M$,
the dual space $M^{*}$ of continuous linear functionals on $M$
can be identified with $M$ through the bicomplex inner product
$\scalarmath{\cdot}{\cdot}$.

We close this section by proving a general version of Schwarz's
inequality in a bicomplex Hilbert space.

\begin{theorem} [Bicomplex Schwarz inequality]
Let $\ket{\psi},\ket{\phi}\in M$.  Then
$$|(\ket{\psi},\ket{\phi})|\leq \sqrt{2}\big{\|}
\ket{\psi} \big{\|} \ \big{\|} \ket{\phi} \big{\|}.$$
\end{theorem}
\begin{proof}
From the complex (in $\mC(\bo)$) Schwarz inequality we have
\begin{equation}
|(\ket{\psi}_\mk,\ket{\phi}_\mk)_{\widehat{k}}|^{2}
\leq\big{|}\ket{\psi}_\mk\big{|}^{2}_{k}
\cdot\big{|}\ket{\phi}_\mk\big{|}^{2}_{k} , \quad
\forall \ket{\psi}_\mk,\ket{\phi}_\mk\in V_k. \notag
\end{equation}
Therefore, if $\ket{\psi},\ket{\phi}\in M$, we obtain
from \eqref{norm7} and \eqref{normk1}
\begin{align*}
|(\ket{\psi},\ket{\phi})|
&= |\eo(\ket{\psi}_\mo,\ket{\phi}_\mo)_{\widehat{1}}
+\et(\ket{\psi}_\mt,\ket{\phi}_\mt)_{\widehat{2}}|\\
&= \frac{1}{\sqrt{2}} \sqrt{|(\ket{\psi}_\mo,\ket{\phi}_\mo)_{\widehat{1}}|^{2}
+|(\ket{\psi}_\mt,\ket{\phi}_\mt)_{\widehat{2}}|^{2}} \\
&\leq
 \frac{1}{\sqrt{2}} \sqrt{\big{|} \ket{\psi}_\mo\big{|}^{2}_{1}
\cdot\big{|} \ket{\phi}_\mo\big{|}^{2}_{1}
+ \big{|} \ket{\psi}_\mt\big{|}^{2}_{2}
\cdot\big{|} \ket{\phi}_\mt\big{|}^{2}_{2}} \\
&\leq \sqrt{2}\big{\|} \ket{\psi} \big{\|}
\; \big{\|} \ket{\phi} \big{\|}.
\end{align*}
\end{proof}

\noindent This result is a direct generalization of the bicomplex
Schwarz inequality obtained in~\cite{Rochon3} and proved here
for an arbitrary bicomplex Hilbert space $M$ without
any extra condition on the bicomplex scalar product.

\subsection{Countable $\T$-Modules}

In this section we investigate more specific $\T$-modules,
namely those that have a countable basis.  Such modules
have another important Hilbert subspace.

\subsubsection{Schauder $\T$-Basis}

\begin{definition}
Let $M$ be a normed $\T$-module. We say that $M$ has a
\textbf{Schauder $\T$-basis} if there exists a countable set
$\{\ket{m_1} \dots \ket{m_l} \dots\}$ of elements of $M$
such that every element $\ket{\psi}\in M$ admits a unique
decomposition as the sum of a convergent series
$\ket{\psi}=\sum_{l=1}^{\infty}w_l\ket{m_l}$, $w_l\in\T$.
\end{definition}

\noindent
If $\{\ket{m_l}\}$ is a Schauder $\T$-basis in $M$, it follows
that $\sum_{l=1}^{\infty}w_l\ket{m_l}=0$
if and only if $w_l=0$, $\forall l\in\mathbb{N}$.
Moreover, if a normed $\T$-module $M$ has a \textbf{Schauder $\T$-basis},
then the vector space $M'$ is automatically a normed vector
space with the following classical Schauder basis:
\begin{equation}
\{\ket{m_1}_\mo, \ket{m_1}_\mt \dots \ket{m_l}_\mo, \ket{m_l}_\mt
\dots\} . \notag
\end{equation}
The normed space $M'$ with a Schauder basis is necessarily of
infinite dimension since it contains subspaces of an arbitrary finite dimension. For more details on Schauder bases see~\cite{Hansen}.

A normed $\T$-module with a Schauder $\T$-basis is called a
\textbf{countable $\T$-module}. For the rest of this section
we will only consider bicomplex Hilbert spaces constructed from
countable $\T$-modules. We now show that in this
context, it is always possible to
construct an orthonormal Schauder $\T$-basis in~$M$.

\begin{theorem}[Orthonormalization]
Let $M$ be a bicomplex Hilbert
space and let $\{ \ket{s_l} \}$ be an arbitrary Schauder $\T$-basis of $M$.
Then $\{ \ket{s_l} \}$ can always be orthonormalized.
\label{Theo3.3.1}
\end{theorem}
\begin{proof}
Using Theorem~\ref{SD} and Corollary~\ref{Hilbert1}, we can write
$M=V_1\oplus V_2$, where $V_1$ and $V_2$ are the associated Hilbert spaces.
Hence $\{ \ket{s_l}_\mk \}$ is a Schauder basis of $V_k$ ($k=1, 2$).
Let $\{\ket{s'_l}_\mk \}$ be the orthonormal basis constructed
from $\{ \ket{s_l}_\mk \}$ in $V_k$~\cite[p.~59]{Hansen}.
For all $l\in\mathbb{N}$ we have that
$$\scalarmath{\ket{s'_l}_\mo + \ket{s'_l}_\mt}
{\ket{s'_l}_\mo + \ket{s'_l}_\mt}
= \eo\big{|} \ket{s'_l}_\mo \big{|}_{1}^2
+ \et\big{|} \ket{s'_l}_\mt \big{|}_{2}^2=1$$
and
$$\scalarmath{\ket{s'_l}_\mo + \ket{s'_l}_\mt}
{\ket{s'_p}_\mo + \ket{s'_p}_\mt}=0$$
if $l\neq p$.  From this we conclude that the set
$\{ \ket{s'_l}_\mo + \ket{s'_l}_\mt \}$
is an orthonormal basis in $M$.
\end{proof}

\medskip\noindent
It is interesting to note that the normalizability of kets requires
that the scalar product belongs to $\D^+$.  To see this, let us write
$\scalarmath{\ket{m_1}}{\ket{m_1}}=a_\hh\ee+a_\hhh\eee$ with
$a_\hh,a_\hhh\in\R$, and let
\begin{equation}
\ket{m'_1}=(z_\hh\ee+z_\hhh\eee)\ket{m_1}, \notag
\end{equation}
with $z_\hh,z_\hhh\in\C(\ii)$ and $z_\hh\neq0\neq z_\hhh$. We get
\begin{align*}
\scalarmath{\ket{m'_1}}{\ket{m'_1}}
&=(|z_\hh|^2\ee+|z_\hhh|^2\eee)\scalarmath{\ket{m_1}}{\ket{m_1}}\\
&=(|z_\hh|^2\ee+|z_\hhh|^2\eee)(a_\hh\ee+a_\hhh\eee)\\
&=c_\hh a_\hh\ee+c_\hhh a_\hhh\eee,
\end{align*}
with $c_\h{k}=|z_\h{k}|^2\in\R^+$. The normalization condition of
$\ket{{m}_1'}$ becomes
\begin{equation}
c_\hh a_\hh\ee+c_\hhh a_\hhh\eee=1, \notag
\end{equation}
or $c_\hh a_\hh=1=c_\hhh a_\hhh$. This is possible only if
$a_\hh > 0$ and $a_\hhh > 0$.  In other words,
$\scalarmath{\ket{m_1}}{\ket{m_1}}\in\D^+$.

\subsubsection{Projection in a Specific Vector Space}

Let $\{\ket{m_1} \dots \ket{m_l} \dots\}$ be a Schauder
$\T$-basis associated with the bicomplex Hilbert space
$\{M, (\cdot,\cdot)\}$.  That is, any element $\ket{\psi}$
of $M$ can be written as
\begin{equation}
\ket{\psi} = \sum_{l=1}^{\infty} w_l\ket{m_l} ,\label{2.22}
\end{equation}
with $w_l \in \T$.  As was shown in~\cite{Rochon3} for the
finite-dimensional case, an important subset $V$ of $M$ is
the set of all kets for which all $w_l$ in \eqref{2.22}
belong to $\C(\ii)$. It is obvious that $V$ is a non-empty
normed vector space over complex numbers with Schauder basis
$\{\ket{m_1} \dots \ket{m_l} \dots\}$. Let us write
$w_l = \ee {z_{l\hh}} + \eee {z_{l\hhh}}$, so that
\begin{equation}
\ket{\psi} = \sum_{l=1}^{\infty} w_l\ket{m_l}
= \sum_{l=1}^{\infty} \ee z_{l\hh} \ket{m_l}
+ \sum_{l=1}^{\infty} \eee z_{l\hhh} \ket{m_l} .
\end{equation}
From Theorem \ref{closed} we know that the two
series on the right-hand side separately converge to elements
of $V_1$ and $V_2$.  However, it is not obvious that
$\sum_{l=1}^{\infty} {z_{l\h{k}}}\ket{m_l}$ ($k=1, 2$)
converges since $\e{k}$ is not invertible. To prove this result,
we need the following theorem.
\begin{theorem}
Let $\{\ket{\psi_n}\}$ be an orthonormal sequence in the bicomplex
Hilbert space $M$ and let $\{\alpha_n\}$ be a sequence of bicomplex numbers.
Then the series $\sum_{n=1}^{\infty} \alpha_n\ket{\psi_n}$ converges in $M$
if and only if $\sum_{n=1}^{\infty} |\alpha_n|^2$ converges in $\mathbb{R}$.
\label{Theo4.6}
\end{theorem}
\begin{proof}
The series $\sum_{n=1}^{\infty} \alpha_n\ket{\psi_n}$ converges
if and only if for $k=1, 2$, the series
$$\sum_{n=1}^{\infty} \e{k}\alpha_n\ket{\psi_n}
=\sum_{n=1}^{\infty} P_k(\alpha_n) \ket{\psi_n}_\mk$$
converges.  However, in the Hilbert space $V_k$, it is
well-known~\cite[p.~59]{Hansen} that
$\sum_{n=1}^{\infty} P_k(\alpha_n) \ket{\psi_n}_\mk$
converges if and only if $\sum_{n=1}^{\infty} |P_k(\alpha_n)|^2$
converges in $\mathbb{R}$.  Since
$$|\alpha_n|^2=\frac{|P_1(\alpha_n)|^2+|P_2(\alpha_n)|^2}{2},$$
we find that the series $\sum_{n=1}^{\infty} \alpha_n\ket{\psi_n}$
converges in $M$ if and only if $\sum_{n=1}^{\infty} |\alpha_n|^2$
converges in $\mathbb{R}$.
\end{proof}

\medskip\noindent
From Theorem~\ref{Theo4.6} we see that if
$\{\ket{m_1} \dots \ket{m_l} \dots\}$ is an \textbf{orthonormal}
Schauder $\T$-basis and
\begin{equation}
\sum_{l=1}^{\infty} (\ee {z_{l\hh}}+\eee {z_{l\hhh}})\ket{m_l}
\notag
\end{equation}
converges in $M$, then the series
\begin{equation}
\sum_{l=1}^{\infty} |\ee {z_{l\hh}}+\eee {z_{l\hhh}}|^2
\notag
\end{equation}
converges in~$\R$. In particular,
$\sum_{l=1}^{\infty} |z_{l\h{k}}|^2$ also converges. Hence
$\sum_{l=1}^{\infty} {z_{l\h{k}}}\ket{m_l}$ converges and this allows
to define projectors $P_1$ and $P_2$ from $M$ to $V$ as
\begin{equation}
\P{k}{\ket{\psi}} := \sum_{l=1}^{\infty} {z_{l\h{k}}}\ket{m_l},
\qquad k=1, 2.
\end{equation}
Therefore, any $\ket{\psi} \in M$ can be decomposed
uniquely as
\begin{equation}
\ket{\psi} = \ee \P{1}{\ket{\psi}}
+ \eee \P{2}{\ket{\psi}}.
\label{2.23}
\end{equation}
As in the finite-dimensional case~\cite{Rochon3},
one can easily show that ket projectors
and idempotent-basis projectors (denoted with the
same symbol) satisfy the following, for $k=1,2$:
\begin{align}
\P{k}{s \ket{\psi} + t \ket{\phi}}
= \P{k}{s} \P{k}{\ket{\psi}}
+ \P{k}{t} \P{k}{\ket{\phi}} .\label{2.24}
\end{align}
It will be useful to rewrite \eqref{2.23} as
\begin{equation}
\ket{\psi} = \ee\ket{\psi}_{\widehat{1}} + \eee\ket{\psi}_{\widehat{2}}
= \ket{\psi}_\mo + \ket{\psi}_\mt ,
\label{2.25}
\end{equation}
where $\ket{\psi}_{\widehat{k}}:= \P{k}{\ket{\psi}}$.
Note that the scalar product in Theorem~\ref{Theo4.1}
can be reinterpreted as
\begin{equation}
(\ket{\psi},\ket{\phi})
=\eo(\ket{\psi}_{\widehat{1}},\ket{\phi}_{\widehat{1}})_{\widehat{1}}
+ \et(\ket{\psi}_{\widehat{2}},\ket{\phi}_{\widehat{2}})_{\widehat{2}} .
\end{equation}
The $\T$-norm defined in \eqref{T-norm} can be written as
\begin{equation}
\big{|}\big{|}\ket{\phi}\big{|}\big{|}
= \frac{1}{\sqrt{2}}\sqrt{\big{|}\ket{\phi}_{\widehat{1}}\big{|}^{2}_{1}
+\big{|}\ket{\phi}_{\widehat{2}}\big{|}^{2}_{2}} \, ,
\label{TV-norm}
\end{equation}
where
\begin{equation}
\big{|} \ket{\phi}_{\widehat{k}} \big{|}_{k}^{2}
= \scalarmath{\ket{\phi}_{\widehat{k}}}
{\ket{\phi}_{\widehat{k}}}_{\widehat{k}}
= \scalarmath{\ket{\phi}_\mk}{\ket{\phi}_\mk}_{\widehat{k}}
= \big{|}\ket{\phi}_\mk \big{|}_{k}^{2} .
\label{lib}
\end{equation}

Just to avoid confusion, we recall that a bold index
$\mathbf{k}$ on a ket (like in $\ket{\psi}_\mk$) means that
the ket belongs to $V_k$.  Kets $\ket{\psi}_{\widehat{1}}$
and $\ket{\psi}_{\widehat{2}}$, on the other hand, belong
to~$V$.  Indices~$\widehat{1}$ and~$\widehat{2}$ on a bicomplex
number denote idempotent projection.  The double bar denotes
the $\T$-norm in~$M$ while the single bar denotes standard
norms in $V_k$ or~$V$.

In Corollary~\ref{Hilbert1} we showed that
$(V_k,\scalarmath{\cdot}{\cdot}_{\widehat{k}})$
is a Hilbert space for $k=1,2$. The same is true with $V$.

\begin{theorem}
Let $M$ be a bicomplex Hilbert space with orthonormal
Schauder $\T$-basis $\{\ket{m_l} \}$.  Then
$(V, \scalarmath{\cdot}{\cdot}_{\widehat{k}})$ is a
complex $($in $\mC(\bo))$ Hilbert space for $k=1,2$.%
\label{VHilbert}
\end{theorem}
\begin{proof}
It is easy to see that $(V, \scalarmath{\cdot}{\cdot}_{\widehat{k}})$
is a normed space over $\mC(\bo)$ for $k=1,2$.
Without lost of generality, let us consider the case $k=1$.
Let $\{\ket{\psi_n}\}\in V$ be a Cauchy sequence with respect to
norm $|\cdot|_{1}$ specified in \eqref{lib}.  We can see
that $\{\ket{\psi_n}_\mo \}$
is also a Cauchy sequence. By Corollary~\ref{Corollary4.1},
$\{\ket{\psi_n}_\mo\}$ converges to a ket
$\ket{\psi}_\mo=\eo\sum_{l=1}^{\infty} w_l\ket{m_l}\in V_1$
with respect to the norm $|\cdot|_{1}$, where $\sum_{l=1}^{\infty} w_l\ket{m_l}\in M$. However,
\begin{equation*}
\eo\sum_{l=1}^{\infty} w_l\ket{m_l}
=\sum_{l=1}^{\infty} \e{1}w_l\ket{m_l}
= \sum_{l=1}^{\infty} \e{1}P_1(w_l)\ket{m_l}.
\end{equation*}
Hence, from the discussion after Theorem~\ref{Theo4.6}, we have that
$$\ket{\psi}_\mo = \e{1}\sum_{l=1}^{\infty} P_1(w_l)\ket{m_l} ,\notag$$
where $\sum_{l=1}^{\infty} P_1(w_l)\ket{m_l}\in V$.
Thus $\ket{\psi_n}_{\bold{1}} = \e{1}\ket{\psi_n}\rightarrow
\e{1}\sum_{l=1}^{\infty} P_1(w_l)\ket{m_l}$ whenever $n\rightarrow\infty$.
Hence by \eqref{lib} we get
$$\ket{\psi_n}\rightarrow \sum_{l=1}^{\infty} P_1(w_l)\ket{m_l}\in V$$
whenever $n\rightarrow\infty$.
\end{proof}

\begin{theorem}
Let $M$ be a bicomplex Hilbert space with $\{\ket{m_l}\}$ an
orthonormal Schauder $\T$-basis of $M$. If the norms
$|\cdot|_{1}$ and $|\cdot|_{2}$ are equivalent on $V$,
then the normed vector space $(V,||\cdot||)$ is closed in $M$.
\label{Vclosed}
\end{theorem}
\begin{proof}
Let $\{\ket{\psi_n}\}\in V$ be a converging sequence to a ket $\ket{\psi}\in M$. Thus,
$$\ket{\psi_n}=\sum_{l=1}^{\infty} z_{nl}\ket{m_l}=\ee\ket{\psi_n}_{\widehat{1}} + \eee\ket{\psi_n}_{\widehat{2}}$$
and
$$\ket{\psi}=\sum_{l=1}^{\infty} w_l\ket{m_l}=\ee\ket{\psi}_{\widehat{1}} + \eee\ket{\psi}_{\widehat{2}}$$
where $z_{nl}\in\mC(\bo)$ and $w_l\in\T$, $\forall n,l\in\mathbb{N}$.
Since $\{\ket{\psi_n}\}\in V$, we have that
$\ket{\psi_n}_{\widehat{1}} =\ket{\psi_n}_{\widehat{2}}$,
$\forall n\in\mathbb{N}$. Using \eqref{TV-norm} for the $\T$-norm,
we find that
\begin{equation*}
\big{|}\ket{\psi_n}_{\widehat{k}}-\ket{\psi}_{\widehat{k}}\big{|}_{k}\rightarrow 0
\end{equation*}
whenever $n\rightarrow\infty$, for $k=1,2$ . Hence, $$\ket{\psi}_{\widehat{1}} =\ket{\psi}_{\widehat{2}}$$
since the norms $|\cdot|_{1}$ and $|\cdot|_{2}$ are equivalent on $V$.
\end{proof}

\smallskip
Since all norms on finite-dimensional vector spaces are equivalent, the normed vector space $(V,||\cdot||)$ is always closed in $M$
whenever $M$ is finite-dimensional. In the infinite-dimensional case, it is possible to find a simple condition to obtain
the closure.

\begin{definition}
Let $\{ \ket{m_l} \}$ be an orthonormal Schauder $\T$-basis of $M$ and let $V$ be
the associated vector space.  We say that a scalar product
is $\mC(\bo)$-closed under $V$ if,
$\forall \ket{\psi},\ket{\phi}\in V$, we have
$(\ket{\psi},\ket{\phi})\in\mC(\bo)$.
\end{definition}

We note that the property of being $\mC(\bo)$-closed is
basis-dependent.  That is, a scalar product may be
$\mC(\bo)$-closed under~$V$ defined through a basis
$\{ \ket{m_l} \}$, but not under~$V'$ defined
through a basis $\{ \ket{s_l} \}$.

\begin{corollary}
Let $M$ be a bicomplex Hilbert space M with $\{\ket{m_l}\}$ an
orthonormal Schauder $\T$-basis of $M$. If the scalar product
is $\mC(\bo)$-closed under $V$ then the inner space $(V,||\cdot||)$ is closed in $M$.
\label{CVclosed}
\end{corollary}
\begin{proof}
Equation \eqref{decompEq} is true whether the bicomplex scalar product
is $\mC(\bo)$-closed under~$V$ or not.  When it is $\mC(\bo)$-closed,
we have for $k = 1, 2$
\begin{equation*}
\scalarmath{\ket{\psi}}{\ket{\phi}}_{\widehat{k}}
= \P{k}{\scalarmath{\ket{\psi}}{\ket{\phi}}}
= \scalarmath{\ket{\psi}}{\ket{\phi}} ,
\qquad \forall \ket{\psi},\ket{\phi}\in V .
\end{equation*}
Hence, $|\cdot|_{1}=|\cdot|_{2}$ and by Theorem \ref{Vclosed}
the inner space $(V,||\cdot||)$ is closed in~$M$.
\end{proof}

\section{The Harmonic Oscillator}

Complex Hilbert spaces are fundamental tools of
quantum mechanics.  We should therefore expect that
bicomplex Hilbert spaces should be relevant to any
attempted generalization of quantum mechanics to
bicomplex numbers.  Steps towards such a generalization
were made in~\cite{GMR, Rochon2,Rochon3}.  In~\cite{GMR},
in particular, the problem of the bicomplex quantum harmonic
oscillator was investigated in detail.

To summarize the main results obtained, let us first
recall the function space introduced in~\cite{GMR}.
Let $n$ be a nonnegative integer and let $\alpha$
be a positive real number.  Consider the following
function of a real variable $x$:
\begin{equation}
f_{n,\alpha} (x) := x^n \, \exp{-\alpha x^2} .
\label{6.1}
\end{equation}
Let $S$ be the set of all finite linear combinations
of functions $f_{n,\alpha} (x)$, with complex
coefficients.  Furthermore, let a bicomplex function
$u(x)$ be defined as
\begin{equation}
u(x) = \ee u_{\hh} (x) + \eee u_{\hhh} (x) ,
\label{6.2}
\end{equation}
where $u_{\hh}$ and $u_{\hhh}$ are both in~$S$.  The set of
all functions $u(x)$ is a $\T$-module, denoted by $M_S$.

Let $u(x)$ and $v(x)$ both belong to $M_S$.
We define a mapping $(u, v)$ of this pair of
functions into~$\D^+$ as follows:
\begin{equation}
\left(u, v \right)
:= \int_{-\infty}^{\infty} u^{\dagger_3}(x) v(x) \rmd x
= \int_{-\infty}^{\infty} \left[ \ee \bar{u}_{\hh} (x)
v_{\hh} (x) + \eee \bar{u}_{\hhh} (x) v_{\hhh} (x) \right] \rmd x .
\label{6.3}
\end{equation}
It is not hard to see that \eqref{6.3}
is always finite and satisfies all the
properties of a bicomplex scalar product.

Let $\xi = \ee \xi_\hh + \eee \xi_\hhh$ be in $\D^+$
and let us define two operators $X$ and $P$
that act on elements of $M_S$ as follows:
\begin{equation}
X \{u(x)\} := x u(x) , \qquad
P \{u(x)\} := -\ii \hbar \xi \frac{\rmd u(x)}{\rmd x} .
\label{6.4}
\end{equation}
It is not difficult to show that
\begin{equation}
[X, P] = \ii \hbar \xi I .
\label{6.5}
\end{equation}
Note that the action of $X$ and $P$ on
elements of $M_S$ always yields elements of $M_S$.
That is, $X$ and $P$ are defined on all of $M_S$.

Let $m$ and $\omega$ be two positive real numbers.
We define the bicomplex harmonic oscillator Hamiltonian
as follows:
\begin{equation}
H := \frac{1}{2m} P^2 + \frac{1}{2} m \omega^2 X^2 .
\label{6.6}
\end{equation}
The problem of the bicomplex quantum harmonic oscillator
consists in finding the eigenvalues and eigenfunctions
of~$H$.

That problem was solved in~\cite{GMR}.  The results can
be summarized as follows.  Let $\theta_k$ ($k = \hh, \hhh$)
be defined as
\begin{equation}
\theta_k := \sqrt{\frac{m \omega}{\hbar \xi_k}} \, x .
\label{6.7}
\end{equation}
Bicomplex harmonic oscillator eigenfunctions
can then be written as (the most general eigenfunction
would have different $l$ indices in the two terms):
\begin{align}
\phi_l (x) &= \ee \phi_{l\hh} + \eee \phi_{l\hhh} \notag\\
&= \ee\left[ \sqrt{\frac{m\omega}{\pi\hbar\xi_{\hh}}}
\frac{1}{2^ll!} \right]^{1/2}
\rme^{-\theta_\hh^2/2} H_l(\theta_\hh)
+ \eee\left[ \sqrt{\frac{m\omega}{\pi\hbar\xi_{\hhh}}}
\frac{1}{2^ll!} \right]^{1/2}
\rme^{-\theta_\hhh^2/2} H_l(\theta_\hhh) ,
\label{6.8}
\end{align}
where $H_l$ are Hermite polynomials \cite{Marchildon}.  Equation \eqref{6.8}
can be written more succinctly as
\begin{equation}
\phi_l (x) = \left[ \sqrt{\frac{m\omega}{\pi\hbar\xi}}
\frac{1}{2^l l!} \right]^{1/2} \rme^{-\theta^2/2}
H_l(\theta) ,
\label{6.9}
\end{equation}
where
\begin{equation}
\theta := \ee \theta_{\hh} + \eee \theta_{\hhh}
\quad \mbox{and} \quad
H_l(\theta) := \ee H_l(\theta_{\hh})
+ \eee H_l(\theta_{\hhh}) .
\label{6.10}
\end{equation}
Finally, it was shown in~\cite{GMR} that
$\tilde{M}$, the collection of all finite linear
combinations of bicomplex functions $\phi_l (x)$,
with bicomplex coefficients, is a $\T$-module.
Specifically,
\begin{equation}
\tilde{M} := \left\{\sum_{l} w_l \phi_l (x)
~|~w_l \in\T\right\} . \label{6.11}
\end{equation}

Since $\tilde{M}$ only involves finite linear
combinations of the functions $\phi_l$, it is not
complete, as was pointed out in~\cite{GMR}.  With the
methods developed in this paper, however, we can
extend $\tilde{M}$ to a complete module, in fact to
a bicomplex Hilbert space.

Just like in Definition~\ref{Definition3.1}, we can define two
vector spaces $\tilde{V}_1$ and $\tilde{V}_2$ as
$\tilde{V}_1 = \ee \tilde{M}$ and
$\tilde{V}_2 = \eee \tilde{M}$.  From \eqref{6.8}
and \eqref{6.11}, it is clear that $\tilde{V}_1$
contains all functions $\ee \phi_{l\hh}$ and $\tilde{V}_2$
contains all $\eee \phi_{l\hhh}$.  Now the functions
$\phi_{l\hh}$ and $\phi_{l\hhh}$ are normalized eigenfunctions of
the usual quantum harmonic oscillator (with $\hbar$
replaced by $\hbar \xi_\hh$ or $\hbar \xi_\hhh$).  It is
well-known~\cite{Szekeres} that, as a Schauder basis,
these eigenfunctions generate $L^2(\R)$.

Let $u(x)$ be defined as in \eqref{6.2}, except that
$u_\hh (x)$ and $u_\hhh (x)$ are both
taken as $L^2(\R)$ functions.
Clearly, the set of all $u(x)$ makes up a $\T$-module,
which we shall denote by $M$.  With the scalar product
\eqref{6.3}, $M$ becomes a bicomplex pre-Hilbert space.
Since $L^2(\R)$ is complete we get from
Theorem~\ref{HilbertT}:

\begin{corollary}
$M$ is a bicomplex Hilbert space.
\end{corollary}

\noindent The following theorem shows that the $\phi_l$
make up a Schauder $\T$-basis of~$M$.
\begin{theorem}
Any bicomplex function $u$ in $M$ can be expanded uniquely as
$u = \sum_{l=0}^\infty w_l \phi_l$, where $w_l$ is a
bicomplex number and $\phi_l$ is given in \eqref{6.8}.
\end{theorem}
\begin{proof}
Let $u = \ee u_\hh + \eee u_\hhh$.  Since $u_\hh$
and $u_\hhh$ both belong to $L^2(\R)$ and since
the functions $\phi_{l\hh}$ and $\phi_{l\hhh}$ are
Schauder bases of $L^2(\R)$, one can write
\begin{equation}
u_\hh = \sum_{l=0}^\infty c_{l\hh} \phi_{l\hh} , \qquad
u_\hhh = \sum_{l=0}^\infty c_{l\hhh} \phi_{l\hhh} , \notag
\end{equation}
where $c_{l\hh}$ and $c_{l\hhh}$ belong to $\mathbb{C}(\ii)$.
Letting $w_l = \ee c_{l\hh} + \eee c_{l\hhh}$, the desired
expansion follows.  It must be unique since otherwise,
either $u_\hh$ or $u_\hhh$ would have two different
expansions.
\end{proof}

\section{Conclusion}

We have derived a number of new results on
infinite-dimensional bicomplex modules and Hilbert spaces,
including a generalization of the Riesz representation
theorem for bicomplex continuous linear functionals and a
general version of the bicomplex Schwarz inequality.
The perspective of further investigating the extent to
which quantum mechanics generalizes to bicomplex numbers
motivates us in developing additional mathematical tools
related to infinite-dimensional bicomplex Hilbert spaces
and operators acting on them.  We believe that results
like the Riesz-Fischer theorem and the spectral theorem
can also be extended to infinite-dimensional Hilbert spaces.

\section*{Acknowledgments}
DR is grateful to the Natural Sciences and
Engineering Research Council of Canada for financial
support.  RGL would like to thank the Qu\'{e}bec
FQRNT Fund for the award of a postgraduate
scholarship.

\bibliographystyle{amsplain}

\end{document}